\documentclass[11pt,a4paper]{article}

\usepackage[T1]{fontenc}
\usepackage[british]{babel}
\usepackage{amsmath,amssymb,amsthm}
\usepackage[margin=1in]{geometry}
\usepackage{microtype}
\usepackage[hidelinks]{hyperref}

\newtheorem{theorem}{Theorem}[section]

\newtheorem{corollary}[theorem]{Corollary}
\newtheorem{lemma}[theorem]{Lemma}
\theoremstyle{definition}

\theoremstyle{remark}

\newcommand{\Z}{\mathbb{Z}}

\title{The dilate rollercoaster}
\author{Luke Pebody}
\date{\today}

\begin{document}

\maketitle

\begin{abstract}
    Given a finite subset $A\subseteq\Z$ and a multiset of non-zero integers
    $\Lambda=\{\lambda_1, \ldots, \lambda_k\}$, the \emph{dilate sum}
    $\Lambda\cdot A$ is defined to be
    \[\Lambda\cdot A
    =\lambda_1\cdot A+\lambda_2\cdot A+\cdots+\lambda_k\cdot A
    =\{\lambda_1a_1+\cdots+\lambda_ka_k: a_1, \ldots, a_k\in A\}.\]

    In a recent paper, Kravitz~\cite{Kravitz} determines for 
    which pairs of multisets $(\Lambda_1, \Lambda_2)$ a finite
    set $A$ exists satisfying $|\Lambda_1\cdot A|>|\Lambda_2\cdot A|$.

    For a non-zero real number $r$, 
    say that $\Lambda_2$ \emph{$r$-dominates} $\Lambda_1$ if we can partition some of the
    elements of $\Lambda_2$ into blocks summing to the elements of 
    $r\cdot\Lambda_1$. Then for all finite subsets $A\subseteq\Z$, 
    $\Lambda_2\cdot A$ contains a translate of $r\cdot(\Lambda_1\cdot A)$, and hence
    $|\Lambda_1\cdot A|\le|\Lambda_2\cdot A|$ for all $A$. 
    Kravitz~\cite{Kravitz} proves the converse: if $\Lambda_2$ does not $r$-dominate
    $\Lambda_1$ for any $r$, then there exists a finite set $A$ satisfying 
    $|\Lambda_1\cdot A|>|\Lambda_2\cdot A|$.

    For example, $\{2, 3, 5, 7\}$ does not contain two disjoint subsets whose sums have
    a ratio of 7, and so there exists a finite $A\subseteq\Z$ for which
    $|A+7\cdot A|>|2\cdot A+3\cdot A+5\cdot A+7\cdot A|$.

    We extend Kravitz's method to apply it to multiple inequalities and
    multiple sets, showing that for multisets $\Lambda_{i,j}$ with $1\le i\le k$ 
    and $1\le j\le\ell$, there exist finite sets $A_1, \ldots, A_\ell\subseteq\Z$ such that
    for all $1\le i<k$
    \[|\Lambda_{i,1}\cdot A_1+\cdots+\Lambda_{i,\ell}\cdot A_\ell|
    >|\Lambda_{i+1,1}\cdot A_1+\cdots+\Lambda_{i+1,\ell}\cdot A_\ell|\]
    provided there do not exist indices $1\le i_1<i_2\le k$ and a non-zero scalar $r$ such that 
    $\Lambda_{i_2,j}$ $r$-dominates $\Lambda_{i_1,j}$ for all $1\le j\le\ell$.

    Along the way, we resolve the analogue of Nathanson's Problem 1
    for his Theorem 4~\cite{NathansonComparison},
    as well as the problems in Section 1.4.3 of 
    Kravitz's paper~\cite{Kravitz}.
\end{abstract}

\section{Introduction}\label{sec:introduction}
For subsets $A, B$ of an abelian group, we denote the sumset
\[A+B:=\{a+b:a\in A,b\in B\}.\]

For a positive integer $h$, we denote the $h$-fold iterated
sumset of $A$ by 
\[hA:=\underbrace{A+\cdots+A}_h=\{a_1+\cdots+a_h:a_1,\ldots,a_h\in A\},\]
and let $0A=\{0\}$.

For a non-zero integer $\lambda$, we denote the product of $\lambda$ and
$A$ by $\lambda\cdot A:=\{\lambda a:a\in A\}$. Given a multiset $\Lambda=\{\lambda_1, \ldots, \lambda_k\}$ of non-zero integers,
we denote the $\Lambda$-dilate of $A$ by 
\[\Lambda\cdot A=\lambda_1\cdot A+\lambda_2\cdot A+\cdots+\lambda_k\cdot A=\{\lambda_1a_1+\cdots+\lambda_ka_k: a_1, \ldots, a_k\in A\},\]
and let $\emptyset\cdot A=\{0\}$.

For a positive integer $n$, we write $[n]:=\{1, \ldots, n\}$ and we denote the cyclic group
of $n$ elements by $\Z_n:=\Z/n\Z$. We will make use of the Chinese Remainder Theorem
particularly in the statement that if $p_1, \ldots, p_k$ are distinct primes with 
product $q$, then  $\Z_q\cong\Z_{p_1}\times\ldots\times\Z_{p_k}$.

Following Kravitz~\cite{Kravitz}, for a fixed non-zero real number $r$, 
we will say that the multiset 
$\Lambda_2=\{\lambda_{2,1}, \ldots, \lambda_{2,\ell}\}$
\emph{$r$-dominates} $\Lambda_1=\{\lambda_{1,1}, \ldots, \lambda_{1,k}\}$
if there are disjoint subsets $S_1, \ldots, S_k$
of $[\ell]$ such that
\[\forall j\in[k]: r\lambda_{1,j}=\sum_{h\in S_j}\lambda_{2,h},\]
and say $\Lambda_2$ \emph{dominates} $\Lambda_1$ if
$\Lambda_2$ $r$-dominates $\Lambda_1$ for some non-zero real number $r$.

For example, $\{2, 3, 5\}$ $1$-dominates $\{2, 8\}$, $-2$-dominates $\{-1, -4\}$
and $\frac12$-dominates $\{4, 10\}$. If $\Lambda_2$ $r$-dominates $\Lambda_1$
in this way, then for all subsets $A$ of $\Z$, 
$\Lambda_2\cdot A$ contains a translate of 
$r\cdot(\Lambda_1\cdot A)$ and hence 
$|\Lambda_1\cdot A|\le|\Lambda_2\cdot A|$. Kravitz~\cite{Kravitz} shows 
that domination is also necessary.

\begin{theorem}[Theorem 1.2 of~\cite{Kravitz}]
    Let $\Lambda_1$ and $\Lambda_2$ be two multisets of non-zero integers. 
    Then the dilate inequality $|\Lambda_1\cdot A|\le|\Lambda_2\cdot A|$
    holds for all finite sets $A\subseteq\Z$ if and only if $\Lambda_2$
    dominates $\Lambda_1$.
\end{theorem}

Kravitz~\cite{Kravitz} also asks questions about multi-legged sumset races. 
In particular, for a fixed positive
integer $H$ he considers the ordered pairs $(s, t)$ with $0\le t\le s\le H-1$.
Suppose that \[(s_1, t_1),\ldots,(s_{H(H+1)/2},t_{H(H+1)/2})\] is an ordering
of these pairs in which $(s,t)$ occurs after both $(s+1, t)$ and (if $s>t$)
$(s, t+1)$. Must there exist a finite set $A\subseteq\Z$ such that
\[|s_1A-t_1A|>|s_2A-t_2A|>\cdots>|s_{H(H+1)/2}A-t_{H(H+1)/2}A|?\]

The theorem above implies that no single one of the inequalities precludes the
existence of such a set $A$. We slightly extend the method in~\cite{Kravitz} and 
show that there must indeed exist such an $A$. We in fact show that a collection 
of dilate sums can be put in any prescribed order, 
subject to domination requirements.

\begin{theorem}\label{t:main}
    Suppose $\Lambda_1, \ldots, \Lambda_k$ are multisets of non-zero integers
    and suppose that for no $i<j$ does $\Lambda_j$ dominate $\Lambda_i$.

    Then there exists a finite set $A\subseteq\Z$ such that 
    \[|\Lambda_1\cdot A|>|\Lambda_2\cdot A|>\cdots>|\Lambda_k\cdot A|.\]
\end{theorem}

Kravitz~\cite{Kravitz} also poses the following two-set version: Suppose that
\[(s_1, t_1), \ldots, (s_{H^2}, t_{H^2})\] is some ordering of the ordered pairs $(s, t)$
with $0\le s, t<H$
in which $(s, t)$ always occurs after $(s+1, t)$ and $(s, t+1)$.
Must there exist finite sets $A, B\subseteq\Z$ such that 
\[|s_1A+t_1B|>|s_2A+t_2B|>\cdots>|s_{H^2}A+t_{H^2}B|?\]

We give a positive answer to this question by extending Theorem~\ref{t:main}
to multiple sets.

\begin{theorem}\label{t:main-multiple-sets}
    Let $\Lambda_{i,j}$ for $1\le i\le k$, $1\le j\le\ell$ be multisets of non-zero integers 
    such that for no
    $1\le i_1<i_2\le k$ does there exist a non-zero real number $r$ such that 
    $\Lambda_{i_2,j}$ $r$-dominates $\Lambda_{i_1,j}$ for every $1\le j\le\ell$.

    Then there exist finite sets $A_1, A_2, \ldots, A_\ell\subseteq\Z$ such that 
    \[
    |\Lambda_{1,1}\cdot A_1+\Lambda_{1,2}\cdot A_2+\cdots+\Lambda_{1,\ell}\cdot A_\ell|>
    \cdots>
    |\Lambda_{k,1}\cdot A_1+\Lambda_{k,2}\cdot A_2+\cdots+\Lambda_{k,\ell}\cdot A_\ell|.
    \]
\end{theorem}

The key ingredient in these results will be the following lemma, which is an extension
of Proposition 3.1 from~\cite{Kravitz}.

\begin{lemma}\label{l:mainstep}
    Suppose $\Lambda_1, \ldots, \Lambda_k$ and $\Lambda^{*}$ are 
    multisets of non-zero integers and that for no $i$ does $\Lambda^{*}$ dominate
    $\Lambda_i$.

    Then for all positive integers $N$ and positive reals $\epsilon$, there exists 
    a positive integer $q$ and a subset $A$ of $\Z_q$ such that:
    \begin{itemize}
        \item $q$ is coprime to $N$,
        \item for all $1\le i\le k$, $\Lambda_i\cdot A=\Z_q$,
        \item $|\Lambda^{*}\cdot A|<\epsilon q$.
    \end{itemize}
\end{lemma}

This result will also allow us to give an exact characterisation of when the conclusion of
Theorem~4 of Nathanson~\cite{NathansonComparison} holds.

\begin{theorem}\label{t:nathanson-generalisation}
    Suppose $\Lambda_1, \ldots, \Lambda_k$ and $\Phi_1, \ldots, \Phi_\ell$ are
    multisets of non-zero integers.

    The following are equivalent:
    \begin{enumerate}
        \item No $\Phi_j$ dominates any $\Lambda_i$ for $1\le i\le k$ and $1\le j\le\ell$.
        \item For every $\epsilon>0$ and $c>0$, there exist an integer $m>c$ 
              and a subset $A$ of $\Z_m$ such that 
              \begin{align*}
                &\Lambda_i\cdot A=\Z_m \quad (1\le i\le k)\textrm{ and}\\
                &|\Phi_j\cdot A|<\epsilon m \quad (1\le j\le\ell).
              \end{align*}
    \end{enumerate}
\end{theorem}

In Section~\ref{s:methodology}, we prove Lemma~\ref{l:mainstep} and 
extend it to Theorem~\ref{t:nathanson-generalisation}. Then in
Section~\ref{s:end-of-proof} we deduce Theorem~\ref{t:main} from 
Lemma~\ref{l:mainstep}, and we extend it to 
Theorem~\ref{t:main-multiple-sets} in Section~\ref{s:multiple-sets}. 
Finally, in Section~\ref{s:examples},
we give some specific examples, solving the aforementioned problems of
Kravitz~\cite{Kravitz}.

\section{Key construction}\label{s:methodology}
We adapt Kravitz's construction~\cite{Kravitz} to prove Lemma~\ref{l:mainstep}.

For this section, fix multisets of non-zero integers 
$\Lambda_1, \ldots, \Lambda_k$ and $\Lambda^{*}$
such that for no $i$ does $\Lambda^{*}$ dominate
$\Lambda_i$. Since every multiset dominates the empty multiset, we know that 
each $\Lambda_i$ is non-empty.

Let the elements of $\Lambda_i$ be $\{\lambda_{i,1}, \ldots, \lambda_{i,k_i}\}$
and the elements of $\Lambda^{*}$ be $\{\lambda^{*}_1, \ldots, \lambda^{*}_m\}$.

Let us denote by $C$ the set of coefficient pairs $(i, j)$ where $1\le i\le k$ and 
$1\le j\le k_i$.

\begin{lemma}\label{l:buildingblock}
    Given any subset $P$ of $[m]$ and any choice of elements $c_h\in C$ 
    for $h\in P$,
    there exist rationals $(q_c: c\in C)$ such that 
    \begin{align*}
        \sum_{j=1}^{k_i}\lambda_{i,j}q_{(i,j)}&=1\quad (1\le i\le k),\\
        \sum_{h\in P}\lambda_h^*q_{c_h}&=0.
    \end{align*}
\end{lemma}

\begin{proof}
    For all $1\le i\le k$ and $1\le j\le k_i$, let us write $\mu_{i,j}$
    for the sum of $\lambda^{*}_h$ over all $h\in P$ with $c_h=(i,j)$. 
    Then $\sum_{h\in P}\lambda^{*}_hq_{c_h}=\sum_{i,j}\mu_{i,j}q_{(i,j)}$.

    Note that by the definition of domination, for no $i$ is the vector
    $(\mu_{i,1}, \ldots, \mu_{i,k_i})$ a non-zero scalar multiple of
    $(\lambda_{i,1}, \ldots, \lambda_{i,k_i})$. Were it so, the disjoint sets
    $\{h\in P:c_h=(i,j)\}, 1\le j\le k_i$ would witness that $\Lambda^{*}$
    dominates $\Lambda_i$.

    Define $(q_c:c\in C)$ by $q_{(i,1)}=\frac{1}{\lambda_{i,1}}$ for all 
    $1\le i\le k$ and $q_{(i,j)}=0$ if $j>1$. Then clearly 
    $\sum_{j=1}^{k_i}\lambda_{i,j}q_{(i,j)}=1$ for all $1\le i\le k$.
    Let us write \[S:=\sum_{i,j}\mu_{i,j}q_{(i,j)}=\sum_i\frac{\mu_{i,1}}{\lambda_{i,1}}.\]

    If $S=0$, we are done. Otherwise, there is some $i$ with $\mu_{i,1}\ne 0$.
    For this $i$, there exists $j$ with $2\le j\le k_i$ such that 
    $\frac{\mu_{i,1}}{\lambda_{i,1}}\ne\frac{\mu_{i,j}}{\lambda_{i,j}}$ and then there
    exist rationals $q, q'$ such that 
    \begin{align*}
        \lambda_{i,1}q+\lambda_{i,j}q'&=0\\
        \mu_{i,1}q+\mu_{i,j}q'&=-S,
    \end{align*}
    and if we replace $q_{(i,1)}$ with $\frac{1}{\lambda_{i,1}}+q$ and 
    $q_{(i,j)}$ with $q'$, we maintain $\sum_{j=1}^{k_i}\lambda_{i,j}q_{(i,j)}=1$,
    but we make $\sum_{i,j}\mu_{i,j}q_{(i,j)}=0$ as required.
\end{proof}

The following corollary records that the equations 
can be solved modulo every sufficiently large prime, uniformly over all choices of 
$P$ and $(c_h)_{h\in P}$. 

\begin{corollary}\label{c:buildingblock2}
    For every sufficiently large prime $p$, every $P\subseteq[m]$ and 
    every choice of $c_h\in C$ for $h\in P$, there exist elements $(q_c:c\in C)$ of $\Z_p$
    such that
    \begin{align*}
        \sum_{j=1}^{k_i}\lambda_{i,j}q_{(i,j)}&=1\quad (1\le i\le k),\\
        \sum_{h\in P}\lambda_h^*q_{c_h}&=0.
    \end{align*}
\end{corollary}

\begin{proof}
    Choose one rational solution for each pair $(P,(c_h:h\in P))$. There are only finitely 
    many denominators among these solutions, so they are invertible modulo every sufficiently
    large prime.
\end{proof}

The following technical lemma provides the inductive step.

\begin{lemma}\label{l:inductive_step}
    Suppose $N$ and $q$ are positive integers with $\gcd(N,q)=1$
    and suppose that functions $(\phi_c:\Z_q\to\Z_q:c\in C)$ satisfy
        \[\sum_{j=1}^{k_i}\lambda_{i,j}\phi_{(i,j)}(x)=x\quad (x\in\Z_q, 1\le i\le k).\]

    Then for every $\epsilon>0$, 
    there exist a positive integer $q'$ coprime to $Nq$
    and functions $(\widetilde\phi_c:\Z_q\times\Z_{q'}\to\Z_q\times\Z_{q'}:c\in C)$
    such that, writing $\pi:\Z_q\times\Z_{q'}\to\Z_q$ for the projection to the first coordinate and
    defining $E$ by 
    \[
E:=
\left\{
\sum_{j=1}^m\lambda_j^*\widetilde\phi_{c_j}(x_j):
\begin{array}{l}
c_1,\ldots,c_m\in C,\quad
x_1,\ldots,x_m\in\Z_q\times\Z_{q'},\\
|\{x_1,\ldots,x_m\}|
=|\{\pi(x_1),\ldots,\pi(x_m)\}|
\end{array}
\right\},
\]
    the following hold:
    \begin{itemize}
        \item $\pi\circ\widetilde\phi_c=\phi_c\circ\pi$ for every $c\in C$,
        \item $\sum_{j=1}^{k_i}\lambda_{i,j}\widetilde\phi_{(i,j)}(x)=x\quad (x\in\Z_q\times\Z_{q'}, 1\le i\le k)$,
        \item $|E|<\epsilon qq'$.
    \end{itemize}
\end{lemma}

\begin{proof}
    Define
    \[T:=\{((x_1, \ldots, x_m), (c_1, \ldots, c_m)):c_1, \ldots, c_m\in C, x_1, \ldots, x_m\in\Z_q\}.\]

    For each $t\in T$, choose a distinct prime $p_t>|T|\epsilon^{-1}$ which does not divide $Nq$ or $\lambda_{i,1}$ for any 
    $i$ and is large enough to satisfy Corollary~\ref{c:buildingblock2}.

    For each such $t=((x_1,\ldots, x_m), (c_1, \ldots, c_m))$, we will define a function
    $\phi_{t,c}:\Z_{q}\times\Z_{p_t}\to\Z_{p_t}$ by defining the behaviour of
    $\phi_{t,c}(u,v)$ depending on $u\in\Z_q$. 
    
    If $u\in\{x_1, \ldots,x_m\}$, set $P_u=\{j\in[m]:x_j=u\}$ and 
    apply Corollary~\ref{c:buildingblock2} with $P_u$ 
    to get elements of $\Z_{p_t}$ 
    $(q_{c,u}:c\in C)$ such that for all $i$
    $\sum_{j=1}^{k_i}\lambda_{i,j}q_{(i,j),u}=1$, and
    $\sum_{h\in P_u}\lambda^{*}_hq_{c_h,u}=0$. Then we define $\phi_{t,c}(u,v)=q_{c,u}v$
    for all $v$.

    If $u\notin\{x_1, \ldots, x_m\}$, define $\phi_{t,(i,1)}(u,v)=(\lambda_{i,1})^{-1}v$
    for each $i$ and $\phi_{t,(i,j)}(u,v)=0$ for $j>1$.

    Then one can check that for each element $(u,v)$ of $\Z_{q}\times\Z_{p_t}$ and 
    for every $i$, \[\sum_{j=1}^{k_i}\lambda_{i,j}\phi_{t,(i,j)}(u,v)=v.\] Suppose
    $y_1, \ldots, y_m$ are elements of $\Z_{p_t}$ such that $y_i=y_j$ whenever
    $x_i=x_j$. Then if, for $u\in\{x_1, \ldots, x_m\}$, we write $v_u$ for the common
    value of the $y_h$ with $h\in P_u$,
    \[
        \sum_{j=1}^m\lambda^{*}_j\phi_{t,c_j}(x_j,y_j)
        =\sum_u\sum_{h\in P_u}\lambda^{*}_h\phi_{t,c_h}(u,v_u)
        =\sum_u\sum_{h\in P_u}\lambda^{*}_hq_{c_h,u}v_u=0.\]

    Let $q'=\prod_{t\in T}p_t$ and 
    $A=\Z_{q}\times\prod_{t\in T}\Z_{p_t}$, and define functions for $c\in C$ by
    \[\widetilde\phi_c(x,(x_t:t\in T))=(\phi_c(x),(\phi_{t,c}(x,x_t):t\in T)).\]

    By the Chinese Remainder Theorem, $\prod_{t\in T}\mathbb Z_{p_t}\cong\mathbb Z_{q'}$,
    so we identify $A$ with $\Z_q\times\Z_{q'}$, while retaining the individual 
    $p_t$ coordinates. Note that none of the $p_t$ divide $Nq$, so $q'$ is coprime to $Nq$.

    By construction $\pi\circ\widetilde\phi_c=\phi_c\circ\pi$ and 
    for every element $a\in A$, 
    $\sum_{j=1}^{k_i}\lambda_{i,j}\widetilde\phi_{(i,j)}(a)=a$ for every $1\le i\le k$. 
    
    Furthermore, for any element $x\in E$, $x$ has an expression of the form 
    $\sum_{j=1}^m\lambda_j^*\widetilde\phi_{c_j}(x_j)$ where $c_1, \ldots, c_m\in C$,
    $x_1,\ldots,x_m\in A$ and $|\{x_1,\ldots,x_m\}|=|\{\pi(x_1),\ldots,\pi(x_m)\}|$.
    
    Let $t=((\pi(x_1), \ldots, \pi(x_m)), (c_1, \ldots, c_m))$. 
    The projection $\pi$ is injective on $\{x_1, \ldots, x_m\}$, so the 
    $p_t$-coordinates of the $x_i$ are constant on each set $P_u$. The vanishing
    identity above therefore shows that the $p_t$-coordinate of $x$ is 0.
    Thus $E$ is contained in the union of the coordinate hyperplanes 
    $\{a\in A:a_t=0\}$ and hence by a simple union bound,
    \[|E|\le\sum_{t}\frac{1}{p_t}|A|<\epsilon |A|=\epsilon qq'.\]
\end{proof}

\begin{proof}[Proof of Lemma~\ref{l:mainstep}]
    If $m=0$, let $q$ be a prime that is larger than $N$, $\frac1\epsilon$ and the absolute
    value of all elements of the $\Lambda_i$ and let $A=\Z_q$. 
    Then for all $i$, $\Lambda_i$ is non-empty (as it is not dominated by $\Lambda^*=\emptyset$),
    and hence $\Lambda_i\cdot A=\Z_q$, while $|\Lambda^{*}\cdot A|=1<\epsilon q$.

    For $m\ge 1$, we will apply Lemma~\ref{l:inductive_step} $m$ times with error $\epsilon/m$.
    Let $q_0=1$ and let $(\phi_{0,c}:\Z_{q_0}\to\Z_{q_0}:c\in C)$ all be the identity map.
    Then \[\sum_{j=1}^{k_i}\lambda_{i,j}\phi_{0,(i,j)}(x)=x \quad (1\le i\le k).\]
    
    Inductively, for $1\le z\le m$, if we have defined $q_{z-1}$ and $(\phi_{z-1,c}:c\in C)$,
    apply Lemma~\ref{l:inductive_step} with $q=q_{z-1}$, with $N$ unchanged and with error
    $\epsilon/m$. Denote the resulting
    modulus by $q'_z$, the new functions by $\phi_{z,c}$ and the exceptional set by 
    $E_z$. Set $q_z=q_{z-1}q'_z$ and identify $\Z_{q_z}$ with $\Z_{q_{z-1}}\times\Z_{q'_z}$.
    Note that inductively each $q_z$ is coprime to $N$.

    Set $q=q_m$ and let $\phi_c=\phi_{m,c}$ for $c\in C$. By construction 
    \[\sum_{j=1}^{k_i}\lambda_{i,j}\phi_{(i,j)}(x)=x \quad (1\le i\le k).\]

    Then if we set $A$ to be the set $\{\phi_c(x):c\in C,x\in\Z_q\}$, we see that for all
    $x\in\Z_q$, $x=\sum_{j=1}^{k_i}\lambda_{i,j}\phi_{(i,j)}(x)\in\Lambda_i\cdot A$,
    so $\Lambda_i\cdot A=\Z_q$ for all $i$.

    Now consider an element $x$ of $\Lambda^*\cdot A$. It has an expression of the form
    $x=\sum_{j=1}^m\lambda^*_j\phi_{c_j}(x_j)$ where $c_j\in C$ and $x_j\in\Z_q$.
    For $1\le j\le m$, let $\pi_j:\Z_{q_j}\to\Z_{q_{j-1}}$ be the projection onto the 
    first coordinate. Then for $0\le j\le m$, define 
    \[\rho_j:=\pi_{j+1}\circ\cdots\circ\pi_m:\Z_q\to\Z_{q_j},\]
    where $\rho_m$ is the identity map.

    Applying these projections to the elements $x_1, \ldots, x_m$, we get a sequence of sets
    \[X_j=\{\rho_j(x_1), \ldots, \rho_j(x_m)\}\] with non-increasing sizes:
    \begin{align*}
    |X_m|\ge|X_{m-1}|\ge\ldots\ge|X_0|.
    \end{align*}
    
    Since $|X_0|=1$ and $|X_m|\le m$, not all of the inequalities $|X_z|\ge|X_{z-1}|$ 
    can be strict, so for some $z$, $|X_z|=|X_{z-1}|$, which
    means the projection from $\Z_{q_z}$ to $\Z_{q_{z-1}}$ is injective on $X_z$ and hence
    \[
        \rho_z(x)=\sum_{j=1}^m\lambda_j^*\phi_{z,c_j}(\rho_z(x_j))\in E_z,
    \]
    where the equality follows from the compatibility of $\phi_{r,c}$ with the projections $\pi_r$.

    Thus $x$ lies in $\rho_z^{-1}(E_z)$, which is of size 
    $|\rho_z^{-1}(E_z)|=|E_z|\frac{q}{q_z}<\frac{\epsilon q}m$. 
    
    Every element of $\Lambda^*\cdot A$ lies in one of these $m$ sets,
    so a simple union bound gives that $|\Lambda^*\cdot A|<\epsilon q$.
\end{proof}

We now deduce Theorem~\ref{t:nathanson-generalisation}.

\begin{proof}[Proof of Theorem~\ref{t:nathanson-generalisation}]
    Suppose $\Lambda_1, \ldots, \Lambda_k$ and $\Phi_1, \ldots, \Phi_\ell$ are
    multisets of non-zero integers.

    If $\Lambda_i$ is empty then every $\Phi_j$ dominates $\Lambda_i$, and taking
    $c>1$ means there are no sets $A$ with $\Lambda_i\cdot A=\Z_m$.

    If $\Phi_j$ dominates non-empty $\Lambda_i$, then its domination scalar is rational; write it as 
    $r=a/b$ for integers $a, b$. Then $b\cdot\Phi_j$ $1$-dominates
    $a\cdot\Lambda_i$, and hence for every finite set $A$, 
    $(b\cdot \Phi_j)\cdot A$ contains a translate of
    $(a\cdot \Lambda_i)\cdot A$.

    If $\Lambda_i\cdot A=\Z_m$, then
    \[
        |\Phi_j\cdot A|
        \ge |b\cdot (\Phi_j\cdot A)|
        \ge |a\Z_m|
        =\frac{m}{\gcd(m,|a|)}\ge\frac{m}{|a|}.
    \]

    Taking $\epsilon<\frac1{|a|}$ contradicts condition 2.

    On the other hand, suppose that no $\Phi_j$ dominates any $\Lambda_i$
    for $1\le i\le k$ and $1\le j\le\ell$, and fix $\epsilon>0$ and $c>0$.
    Let $m_0$ be any prime greater than $c$ and greater than the absolute value of every
    element of each $\Lambda_i$ and set $A_0=\Z_{m_0}$. Since every $\Lambda_i$ is non-empty
    and every one of its coefficients is invertible modulo $m_0$, we have $\Lambda_i\cdot A_0=\Z_{m_0}$.

    For $1\le t\le\ell$, suppose inductively that we have found an integer $m_{t-1}>c$ and a subset $A_{t-1}$ of $\Z_{m_{t-1}}$ such 
    that
    \begin{align*}
    &\Lambda_i\cdot A_{t-1}=\Z_{m_{t-1}} \quad (1\le i\le k)\textrm{ and}\\
    &|\Phi_j\cdot A_{t-1}|<\epsilon m_{t-1} \quad (1\le j<t).
    \end{align*}

    Then apply Lemma~\ref{l:mainstep} with $\Lambda^*=\Phi_t$, with error $\epsilon$ and with $N=m_{t-1}$ 
    to get a positive integer $q$ and a 
    subset $A$ of $\Z_q$ such that:
    \begin{itemize}
        \item $q$ is coprime to $m_{t-1}$,
        \item for all $1\le i\le k$, $\Lambda_i\cdot A=\Z_q$,
        \item $|\Phi_t\cdot A|<\epsilon q$.
    \end{itemize}
    Set $m_t=qm_{t-1}>c$ and $A_t$ to be the image of $A_{t-1}\times A$ under the Chinese
    Remainder Theorem isomorphism $\Z_{m_{t-1}}\times\Z_{q}\cong\Z_{m_t}$.

    For any multiset $\Lambda$, $\Lambda\cdot A_t\cong(\Lambda\cdot A_{t-1})\times(\Lambda\cdot A)$,
    $|\Lambda\cdot A|\le q$ and $|\Lambda\cdot A_{t-1}|\le m_{t-1}$, so it 
    follows that 
    \begin{align*}
    &\Lambda_i\cdot A_t=\Z_{m_t} \quad (1\le i\le k)\textrm{ and}\\
    &|\Phi_j\cdot A_t|=|\Phi_j\cdot A_{t-1}||\Phi_j\cdot A|<(\epsilon m_{t-1})q=\epsilon m_t \quad (1\le j<t)\textrm{ and}\\
    &|\Phi_t\cdot A_t|=|\Phi_t\cdot A_{t-1}||\Phi_t\cdot A|<m_{t-1}(\epsilon q)=\epsilon m_t
    \end{align*}
    so the induction is complete.

    Set $m=m_\ell$ and $A=A_\ell$. Then $m>c$ and $A\subseteq\Z_m$ satisfies
    \begin{align*}
    &\Lambda_i\cdot A=\Z_m \quad (1\le i\le k)\textrm{ and}\\
    &|\Phi_j\cdot A|<\epsilon m \quad (1\le j\le\ell).
    \end{align*}    
\end{proof}

\section{From finite cyclic groups to the integers}\label{s:end-of-proof}
The proof of the following corollary follows a similar line to
Theorem~\ref{t:nathanson-generalisation}, but combines
invocations of Lemma~\ref{l:mainstep} in a slightly different
way.

\begin{corollary}\label{C:to_sequences}
    Suppose $\Lambda_1, \ldots, \Lambda_k$ are multisets of non-zero integers
    and suppose that for no $i<j$ does $\Lambda_j$ dominate $\Lambda_i$.

    Then for all positive reals $\epsilon$, there exists a positive integer $q$
    and a subset $A$ of $\Z_q$ such that for all $1\le i<k$, 
    $|\Lambda_{i+1}\cdot A|<\epsilon|\Lambda_i\cdot A|$.
\end{corollary}

\begin{proof}
    If $k=1$, set $q=1$ and $A=\Z_1$ and both conditions are vacuous. Suppose $k\ge 2$.

    By Lemma~\ref{l:mainstep} with $N=1$ and error $\epsilon$
    applied to $\Lambda_1$ and $\Lambda^{*}=\Lambda_2$, there exists
    a positive integer $q_1$ and a subset $A_1$ of $\Z_{q_1}$ such that
    $\Lambda_1\cdot A_1=\Z_{q_1}$ and $|\Lambda_2\cdot A_1|<\epsilon q_1$.

    Suppose inductively for $2\le t<k$ we have found a positive integer
    $q_{t-1}$ and a subset $A_{t-1}$ of 
    $\Z_{q_{t-1}}$ such that for all $1\le j<t$, 
    \[|\Lambda_{j+1}\cdot A_{t-1}|<\epsilon |\Lambda_j\cdot A_{t-1}|.\]

    Choose $\epsilon_t$ such that 
    $\epsilon_t|\Lambda_{t+1}\cdot A_{t-1}|<\epsilon|\Lambda_t\cdot A_{t-1}|$.
    Apply Lemma~\ref{l:mainstep} with targets $\Lambda_1, \ldots, \Lambda_{t}$,
    exceptional multiset $\Lambda_{t+1}$, $N=q_{t-1}$ and error $\epsilon_t$.
    Then there exists a positive integer $q$
    coprime to $q_{t-1}$ and a subset $A$ of $\Z_{q}$ such that
    $\Lambda_j\cdot A=\Z_{q}$ for all $1\le j\le t$
    and $|\Lambda_{t+1}\cdot A|<\epsilon_tq$.

    Then set $q_{t}=qq_{t-1}$ and let $A_t$ 
    be the image of $A\times A_{t-1}$ under the Chinese Remainder Theorem isomorphism $\Z_q\times\Z_{q_{t-1}}\cong\Z_{q_t}$. 
    Then for all $1\le j<t$,
    \[|\Lambda_{j+1}\cdot A_t|=|\Lambda_{j+1}\cdot A_{t-1}||\Lambda_{j+1}\cdot A|
        =q|\Lambda_{j+1}\cdot A_{t-1}|<q\epsilon |\Lambda_j\cdot A_{t-1}|=\epsilon|\Lambda_j\cdot A_t|.\]
    
    Further, 
    \[
        |\Lambda_{t+1}\cdot A_t|
        =|\Lambda_{t+1}\cdot A_{t-1}||\Lambda_{t+1}\cdot A|
        <q\epsilon_t|\Lambda_{t+1}\cdot A_{t-1}|
        <q\epsilon|\Lambda_t\cdot A_{t-1}|
        =\epsilon|\Lambda_t\cdot A_t|.
    \]

    Taking $q=q_{k-1}$ and $A=A_{k-1}$ completes the proof.
\end{proof}

Now we can deduce Theorem~\ref{t:main}.

\begin{proof}[Proof of Theorem~\ref{t:main}]
    Apply Corollary~\ref{C:to_sequences} to the multisets $\Lambda_1, \ldots, \Lambda_k$ with an
    error value $\epsilon$ we will specify later to get a subset $A$ of $\Z_q$ satisfying
    the required inequalities.

    Let $A^*$ be the subset of the standard representatives $\{0, 1, \ldots, q-1\}$ 
    which is the preimage of $A$ under the canonical projection map. 

    Then it is clear that for every multiset $\Lambda$ of non-zero integers,
    $\Lambda\cdot A^*$ contains at least one preimage of each element of $\Lambda\cdot A$,
    so $|\Lambda\cdot A^*|\ge|\Lambda\cdot A|$.

    Furthermore, if $\Lambda=\{\lambda_1, \lambda_2, \ldots, \lambda_h\}$ and we write 
    $\mathrm{wt}(\Lambda)=|\lambda_1|+\cdots+|\lambda_h|$, then one can show that 
    $(2\mathrm{wt}(\Lambda)+1)|\Lambda\cdot A|\ge |\Lambda\cdot A^*|$. If $\Lambda$ is empty,
    both $|\Lambda\cdot A|$ and $|\Lambda\cdot A^*|$ are equal to 1. Otherwise,
    it is clear that for any element
    $a\in\Lambda\cdot A^*$, $|a|<\mathrm{wt}(\Lambda)q$ and hence there are at most $2\mathrm{wt}(\Lambda)$
    preimages of any given element of $\Lambda\cdot A$.

    Choose the error value $\epsilon$ above, so that it is smaller than $\frac{1}{2\mathrm{wt}(\Lambda_i)+1}$
    for all $1\le i\le k$, 
    \[
        |\Lambda_{i+1}\cdot A^*|
        \le(2\mathrm{wt}(\Lambda_{i+1})+1)|\Lambda_{i+1}\cdot A|
        <(2\mathrm{wt}(\Lambda_{i+1})+1)\epsilon|\Lambda_i\cdot A|
        <|\Lambda_i\cdot A|
        \le |\Lambda_i\cdot A^*|
    \]
    as is required.
\end{proof}

\section{Multiple sets}\label{s:multiple-sets}
In this section, we will use the multiset union operator $\uplus$. 
Given multisets $\Lambda_1$ and $\Lambda_2$, we denote the product
\[\Lambda_1\times\Lambda_2=\biguplus_{j\in\Lambda_1}\left(j\cdot\Lambda_2\right).\]

For example, $\{1,1,2\}\times\{2,4\}=\{2,2,4,4,4,8\}$. 

\begin{lemma}\label{l:multiple-sets}
Suppose $\Lambda_{i,j}$ for $1\le i\le 2$ and $1\le j\le k$ are multisets of at most $t$ non-zero
integers, each of which has absolute value at most $p$.

Fix an integer $N\ge 2tp+1$ and for $1\le i\le 2$ define the multisets
\[
\Phi_i:=
\biguplus_{j=1}^k
\left(\Lambda_{i,j}\times
\{N^{2kj},N^{(2k+1)j}\}\right).
\]

For every non-zero real number $r$, the following are equivalent:
\begin{enumerate}
    \item $\Phi_1$ $r$-dominates $\Phi_2$
    \item $\Lambda_{1,j}$ $r$-dominates $\Lambda_{2,j}$ for all $j$.
\end{enumerate}
\end{lemma}

\begin{proof}
    Given a multiset $S$, write $\sigma(S):=\sum_{s\in S}s$.

    First note that every integer $z$ has at most one expression of the form
    $z=\sum_{i=0}^\ell a_iN^i$ with integer coefficients of absolute value at 
    most $tp$ and with non-zero leading coefficient: suppose we had two
    distinct such expressions $z=\sum_{i=0}^{\ell_1} a_iN^i=\sum_{i=0}^{\ell_2}b_iN^i$
    and let $i$ be the smallest index where $a_i\ne b_i$ (padding the shorter sequence with 
    zeroes if necessary), 
    then $N$ is a divisor of $a_i-b_i$, but $0<|a_i-b_i|\le 2tp<N$, which is impossible.
    
    Now suppose we have submultisets $S, T$ of $\Phi_1$ with $\sigma(S)=N^h\sigma(T)\ne 0$ for
    some $1\le h\le k$. Let us write 
    \begin{align*}
        S&=\biguplus_{j=1}^k\left(N^{2kj}\cdot S_{j,1}\uplus N^{(2k+1)j}\cdot S_{j,2}\right)\\
        T&=\biguplus_{j=1}^k\left(N^{2kj}\cdot T_{j,1}\uplus N^{(2k+1)j}\cdot T_{j,2}\right)
    \end{align*}
    where $S_{j,v}, T_{j,v}\subseteq\Lambda_{1,j}$. We can therefore write down two such expressions 
    for $\sigma(S)=N^h\sigma(T)$, namely:
    \begin{align*}
        \sigma(S)&=\sum_{j=1}^{k}(N^{2kj}\sigma(S_{j,1})+N^{(2k+1)j}\sigma(S_{j,2}))\\
        N^h\sigma(T)&=\sum_{j=1}^{k}(N^{2kj+h}\sigma(T_{j,1})+N^{(2k+1)j+h}\sigma(T_{j,2})).
    \end{align*}

    These two expressions must be the same. Because $1\le h\le k$, the only exponent that appears in
    both expressions is $(2k+1)h$. Indeed, writing $E=\{2k,(2k+1)\}\times[k]$, one checks that
    $E\cap(E+h)=\{(2k+1)h\}$ for $1\le h\le k$. Thus, after discarding any zero-sum layers, 
    $S=N^{(2k+1)h}\cdot S_{h,2}$ and $T=N^{2kh}\cdot T_{h,1}$ for some submultisets $S_{h,2}, T_{h,1}$
    of $\Lambda_{1,h}$ with the same sum.

    Now because $\Phi_1$ $r$-dominates $\Phi_2$, there is a partition of some of 
    the elements of $\Phi_1$ into non-empty disjoint submultisets, such that for each element $\lambda$ of
    $\Phi_2$ there is a corresponding part whose elements sum to $r\lambda$.

    Consider an element $\lambda$ of $\Lambda_{2,j}$. There are corresponding elements
    $\lambda N^{2kj}$ and $\lambda N^{(2k+1)j}$ of $\Phi_2$. Thus there are corresponding parts
    of the domination partition $T$ and $S$, summing to $r\lambda N^{2kj}$ and $r\lambda N^{(2k+1)j}$ 
    respectively. 

    After applying the preceding observation with $h=j$, the surviving low-layer portion has the form
    $N^{2kj}T_{j,1}$ with $\sigma(T_{j,1})=r\lambda$. These $T_{j,1}$ remain disjoint as the occurrences
    of $\lambda$ vary, because the original domination blocks were disjoint. Hence they witness that
    $\Lambda_{1,j}$ $r$-dominates $\Lambda_{2,j}$.

    Conversely, if every $\Lambda_{1,j}$ $r$-dominates $\Lambda_{2,j}$, scale the witnessing blocks by
    $N^{2kj}$ and $N^{(2k+1)j}$. The resulting blocks are disjoint across layers and witness that
    $\Phi_1$ $r$-dominates $\Phi_2$.
\end{proof}

This allows us to determine exactly which inequalities of dilate sums of multiple sets can hold.

\begin{proof}[Proof of Theorem~\ref{t:main-multiple-sets}]
    Choose positive integers $t$ and $p$ such that 
    every $\Lambda_{i,j}$ has at most $t$ elements, each of absolute value
    at most $p$. Then if we set $N=2tp+1$, Lemma~\ref{l:multiple-sets} tells us
    that if we define
    \[\Phi_i=
    \biguplus_{j=1}^{\ell}
\left(\Lambda_{i,j}\times
\{N^{2\ell j},N^{(2\ell+1)j}\}\right)\] 
    for $1\le i\le k$, then the condition in the statement of the theorem implies that
    for no $i_1<i_2$ does $\Phi_{i_2}$ dominate $\Phi_{i_1}$.

    Then Theorem~\ref{t:main} implies the existence of a finite set $A$ such that
    $|\Phi_i\cdot A|>|\Phi_{i+1}\cdot A|$ for all $1\le i<k$. For $1\le j\le\ell$, set
    \[A_j=\{N^{2\ell j}, N^{(2\ell+1)j}\}\cdot A.\]
    Then, for every $1\le i\le k$,
    \[\Phi_i\cdot A=\Lambda_{i,1}\cdot A_1+\Lambda_{i,2}\cdot A_2+\cdots+\Lambda_{i,\ell}\cdot A_\ell.\]

    The required inequalities for the $\Lambda_{i,j}$ now follow from the inequalities for the
    $\Phi_i$.
\end{proof}

\section{Examples}\label{s:examples}
We show here some examples of results that follow from Theorem~\ref{t:main}.
The first two give positive answers to questions which were asked in~\cite{Kravitz}. For non-negative
integer $h$ and multiset $\Lambda$, let us write 
\[h\odot\Lambda:=\biguplus_{u=1}^h\Lambda,\] where the union is empty when $h=0$.

In particular, $h\odot\{1\}$ is the multiset containing $h$ copies of 1. Note that 
$h_1\odot\{1\}$ dominates $h_2\odot\{1\}$ if and only if $h_1\ge h_2$. Note also that for
any finite set $A$, $hA=(h\odot\{1\})\cdot A$.

\begin{corollary}\label{C:sums_of_two_sets}
    Suppose that $(s_1, t_1), \ldots, (s_{H^2}, t_{H^2})$ is an ordering of
    the ordered pairs $(s, t)$ with $0\le s, t\le H-1$ in which $(s, t)$ always occurs after $(s+1, t)$ 
    (assuming $s<H-1$) and $(s, t + 1)$ (assuming $t<H-1$).

    There exist finite sets $A, B\subseteq\Z$ such that 
    \[|s_1A+t_1B|>|s_2A+t_2B|>\ldots>|s_{H^2}A+t_{H^2}B|.\]
\end{corollary}

\begin{proof}
    For $1\le i\le H^2$, let $\Lambda_{i,1}=s_i\odot\{1\}$ and 
    $\Lambda_{i,2}=t_i\odot\{1\}$. Then by the definition of the ordering $(s_i, t_i)$, 
    if $i<j$, at least one of $s_j<s_i$ or $t_j<t_i$ must hold.
    Hence for some $v\in\{1,2\}$, $\Lambda_{j,v}$ does not dominate $\Lambda_{i,v}$.
    
    Thus by Theorem~\ref{t:main-multiple-sets}, there exist finite sets $A_1, A_2$
    such that for all $1\le i<H^2$,
    $|(s_i\odot\{1\})\cdot A_1+(t_i\odot\{1\})\cdot A_2|>|(s_{i+1}\odot\{1\})\cdot A_1+(t_{i+1}\odot\{1\})\cdot A_2|$.

    Setting $A=A_1$ and $B=A_2$, this can be rewritten as 
    \[|s_1A+t_1B|>|s_2A+t_2B|>\ldots>|s_{H^2}A+t_{H^2}B|.\]
\end{proof}

\begin{corollary}\label{C:difference_of_one_set}
    Suppose that $(s_1, t_1), \ldots, (s_{H(H+1)/2}, t_{H(H+1)/2})$ is an ordering of
    the ordered pairs $(s, t)$ with $0\le t\le s\le H-1$ in which $(s,t)$ occurs after $(s+1,t)$ (if $s<H-1$) and
    after $(s,t+1)$ (if $s>t$). 

    There exists a finite set $A\subseteq\Z$ such that
    \[|s_1A-t_1A|>|s_2A-t_2A|>\ldots>|s_{H(H+1)/2}A-t_{H(H+1)/2}A|.\]
\end{corollary}

\begin{proof}
    For $1\le i\le H(H+1)/2$, let $\Lambda_i=(s_i\odot\{1\})\uplus(t_i\odot\{-1\})$.

    If $\Lambda_i$ $r$-dominates $\Lambda_j$ for some non-zero real number $r$ then
     there exists a way of
    partitioning some of the elements of $\Lambda_i$ into parts that sum to
    $r$ times the elements of $\Lambda_j$. We will show that $i\le j$.
    
    If $\Lambda_j$ is empty, then $j=H(H+1)/2$ and hence $i\le j$. We may therefore
    suppose that $\Lambda_j$ is non-empty. Then the 
    parts sum to either $r$ or $-r$ and hence $r$ is a non-zero integer.

    Counting the copies of 1 and -1 in the blocks, if $r>0$, then 
    $s_i\ge rs_j\ge s_j$ and $t_i\ge rt_j\ge t_j$, so $i\le j$. If $r<0$,
    then $s_i\ge |r|t_j\ge t_j$ and $t_i\ge |r|s_j\ge s_j$. Since $s_i\ge t_i$ and
    $s_j\ge t_j$, this again gives $s_i\ge t_i\ge s_j$ and $t_i\ge s_j\ge t_j$, so $i\le j$.

    Thus no $\Lambda_i$ dominates $\Lambda_j$ with $j<i$ and Theorem~\ref{t:main}
    gives the existence of such an $A$.
\end{proof}

\begin{corollary}\label{C:dilates_of_one_set}
    Suppose that $(s_1, t_1), \ldots, (s_n, t_n)$ are integer vectors with 
    non-zero coordinates such that for all $i\ne j$, $(s_i, t_i)$ is not
    parallel to either $(s_j, t_j)$ or $(t_j, s_j)$.

    Then there exists a finite set $A\subseteq\Z$ such that
    \[|s_1\cdot A+t_1\cdot A|>|s_2\cdot A+t_2\cdot A|>\ldots>|s_n\cdot A+t_n\cdot A|.\]
\end{corollary}

\begin{proof}
    For $1\le i\le n$, let $\Lambda_i=\{s_i, t_i\}$.

    If $\Lambda_i$ $r$-dominates $\Lambda_j$, 
    there is a way of partitioning some of the elements of $\Lambda_i$ into 
    blocks that sum to $r$ times the elements of $\Lambda_j$. Since
    each block will need a non-zero sum, both blocks must
    be non-empty. Then, since $\Lambda_i$ and $\Lambda_j$ are the same size,
    each block is a singleton.

    So either $(s_i,t_i)=r(s_j,t_j)$ or $(s_i,t_i)=r(t_j, s_j)$.
    Since $(s_i, t_i)$ is not parallel to either $(s_j, t_j)$ or $(t_j, s_j)$ 
    when
    $i\ne j$, it follows that $i=j$.

    Thus no $\Lambda_i$ dominates any $\Lambda_j$ with $j<i$ and Theorem~\ref{t:main}
    gives the existence of such an $A$.
\end{proof}
\section*{AI}
The author used ChatGPT for a literature search and Codex for proofreading.

\end{document}